Maciej Paszynski
AGH University of Science and Technology
Faculty of Computer Science, Electronics and Telecommunication
Department of Computer Science
al. A. Mickiewicza 30
30-059 Krakow, Poland
tel. +48123283400, fax. +48126175172
e-mail: maciej.paszynski@agh.edu.pl

David Pardo
Department of Applied Mathematics, Statistics and Operational Research, University of the Basque Country UPV/EHU, Leioa, Spain, and
IKERBASQUE (Basque Foundation of Sciences), Bilbao, Spain

Victor Calo
King Abdullah University of Science and Technology, Thuwal, Saudi Arabia

**A direct solver with reutilization of previously-computed LU factorizations for h-adaptive finite element grids with point singularities**

**Abstract**
This paper describes a direct solver algorithm for a sequence of finite element meshes that are *h*-refined towards one or several point singularities. For such a sequence of grids, the solver delivers linear computational cost $O(N)$ in terms of CPU time and memory with respect to the number of unknowns $N$. The linear computational cost is achieved by utilizing the recursive structure provided by the sequence of *h*-adaptive grids with a special construction of the elimination tree that allows for reutilization of previously computed partial LU factorizations over the entire unrefined part of the computational mesh. The reutilization technique reduces the computational cost of the entire sequence of *h*-refined grids from $O(N^2)$ down to $O(N)$. Theoretical estimates are illustrated with numerical results on two- and three-dimensional model problems exhibiting one or several point singularities.

**1. Introduction**
Direct solvers are at the core of many engineering analysis based on the Finite Element Method (FEM) [20]. The finite element (FE) solution process starts by describing a physical phenomenon in terms of a Partial Differential Equation (PDE) with the corresponding boundary and initial conditions over a prescribed domain. Then, this domain is discretized in terms of a finite element mesh, over which the FE solution of the PDE system is obtained by solving a system of linear equations. Finally, the error of the solution is estimated, and if it is above a prescribed tolerance error, the mesh is further refined and the corresponding systems of equations are solved until the solution accuracy is within certain error bounds.
The multi-frontal solver is the state-of-the-art algorithm for solving linear systems of equations [10, 11] when using a direct solver. It is a generalization of the frontal solver algorithm proposed in [9]. The multi-frontal algorithm constructs an assembly tree based on the analysis of the connectivity data. FE are joint into pairs, and fully assembled unknowns are eliminated within frontal matrices associated to multiple branches of the tree. The process is repeated until the root of the assembly tree is reached. Finally, the common interface problem is solved and partial backward substitutions are recursively executed over the assembly tree.
The direct solver algorithm can be generalized to the use of matrix blocks associated with nodes of the computational mesh (called supernodes) [1]. This allows for a reduction of the computational cost related to construction of the elimination tree, since the analysis can then be performed at the level of matrix blocks rather than at the level of particular scalar values. There also exist different implementations of the multi-frontal solver algorithm that target specific computer architectures (see, *e.g.,* [2,3,4,5]). Other significant advances in the area of multi-frontal solvers include the design of a hybrid solver, where the elimination tree is cut at some level, and the remaining Schur complements are submitted to an iterative solver [6]. There also exists a linear computational cost direct solver based on the use of H-matrices [15] with compressed non-diagonal blocks. The main limitation of these solvers is that they produce a non-exact solution.

In this paper we focus on two and three dimensional finite element problems exhibiting point singularities. The corresponding grids are subject to a sequence of *h* refinements towards the existing singularities. The resulting recursive structure of these grids enables one to build a direct solver algorithm that has linear computational cost for a *single* grid from the sequence of *h*-refined grids [19]. In this paper, we describe a special way of constructing the elimination tree that results in a linear computational cost of the direct solver algorithm for the *entire* sequence of *h*-refined grids.

Most available multi-frontal solvers construct the elimination tree by analyzing the structure of the global matrix (e.g., MUMPS solver [12, 13, 14]). In this paper we present an alternative approach, where the elimination tree is constructed based on the direct analysis of the computational mesh. Thus, the input data used during the analysis phase is just the computational mesh. This approach was already proposed in [16], where the elimination tree was constructed by browsing the refinement trees from the smallest to the largest elements following the natural ordering provided by the FE method. Unfortunately, such an ordering results in a non-linear computational cost for *h*-adaptive grids. In our approach, we browse the refinement trees from the largest to the smallest element and we merge frontal matrices from refinement trees adjacent to a common point singularity. By doing so, we achieve linear computational cost for each *h*-refined grid towards one or several point singularities. We have already proposed this ordering in terms of the graph grammar model of the computational mesh [22, 23], however we didn't realize that this ordering implies linear computational cost. Moreover, the elimination tree can be constructed in such a way that previously computed LU factorizations can be reutilized over all unrefined parts the mesh. This feature provides us with a direct solver algorithm that delivers linear computational cost with respect to the number of unknowns not only for a single *h*-refined computational mesh, but also for the entire sequence of meshes constructed by executing several *h*-refinements towards one or many point singularities. Thus, the reutilization technique reduces the computational cost from $O(N^2)$ down to $O(N)$.

Theoretical estimates are illustrated with several numerical results, including 2D problems discretized with the so-called ´radical´ meshes [17, 18] towards one or two singularities, and a 3D Fichera model problem [8].

**2. Elimination Tree of the Direct Solver**

For illustration purposes, we consider the grid described in Figure 1, which describes an initial two elements grid that has been *h*-refined towards a point singularity located in the middle of the bottom edge. We assume a global uniform order of approximation *p* over the entire mesh.

In order to build an elimination tree of the above model problem based on the connectivity information that exists within every FE code, it may seem attractive the idea of building an elimination tree that follows the natural ordering of elements provided by the FE code. Unfortunately, such an approach produces the following undesired situation. In the last step of the solver, all unknowns associated to the vertical mid-edges of the grid (lying at line z=0 in Figure 1) remain untouched, since none of them has been eliminated. As a result, the CPU time cost associated to the elimination performed on this last step of the solver scales as $O(N^{1.5})$, which is prohibitively expensive.

Thus, we propose an alternative approach for building the elimination tree. The solver algorithm browses the refinement trees from the top level down to the leaves, and it uses only one frontal matrix. It first identifies fully assembled nodes located within each level of the refinement elimination tree, eliminates them, and then it iterates the process by going down to the next level.

This elimination tree ensures that the size of a single frontal matrix involved in the solver algorithm remains constant. In the first step, the largest-size elements (namely, elements 1,2,3,4,5, and 6) are eliminated with respect to all interface unknowns (grey stars in Figure 2). After these elements have been eliminated, we obtain a grid with the same topological structure. In the second step of the algorithm, we eliminate elements 7, 8, 9, 10, 11, and 12 with respect to the interface unknowns (dots in Figure 2). This second step has the same cost as the first step, since the number of unknowns that are eliminated is the same in both cases, and the number of interface unknowns is also the same. Then, we iterate the process. We observe that all steps require the same computational cost with the exception of the final step (where we eliminate the remaining elements), which has an even smaller cost than each of the previous steps. Since the cost at each step (level of the elimination tree) is constant, the total cost of the algorithm is proportional to the number of levels, which by grid construction is proportional to the number of unknowns. As a result, we obtain a solver algorithm with linear computational cost with respect to the number of unknowns. For a detailed analysis of the computational cost of an analogous version of the algorithm that browses elements level by level in a reverse bottom-up fashion, we refer to [19].

The main advantage of constructing a top-to-bottom elimination tree (as opposed to the bottom-to-top elimination tree proposed in [19]) is related to the re-utilization of previously computed LU factorizations. If one solves the problem in a given grid and then performs an additional *h*-refinement towards the singular point, the top-to-bottom approach enables a full re-utilization of previously computed LU factorizations, since new refined elements appear at the top of the elimination tree, compare Figure 3. Thus, one only needs to re-compute the LU factorization corresponding to the newly-generated elements. Notice that this is not possible in the bottom-to-top approach proposed in [19], since newly-generated elements appear at the leaves of the elimination tree, which requires re-computation of almost all previously computed partial LU factorizations.

### 3. Theoretical estimates on the computational cost

This section compares the number of floating point operations (NFLOPS) required to solve the entire sequence of *h*-refined grids performed with and without reutilization of previously computed LU factorizations. Again, we consider the radical grid described in Figure 1. We assume a starting grid with eight elements that has been generated by globally refining a two-elements grid. . The sequence of grids is generated by isotropically *h*-refining the two elements closest to the singular point, and iterating the process.

We denote by $L$ to the total number of grids in the sequence. Notice that $L$ is also the number of refinement levels in the last grid from the sequence. Let $T_l$ denotes the execution time of the direct solver algorithm over the *l-th* grid from the sequence of grids. Let $M_l$ denote the memory usage of the direct solver algorithm over the *l-th*. Let $N_l$ denote the number of unknowns at grid number $l$ from the sequence. Let $N = N_L$ be the total number of unknowns at the last *L*-th grid.

As mentioned before (and also proved in [19]), the execution time of the solution of a single grid is linear, which can be expressed as $T_l = c_1 + N_l c_2$ for the *l*-th grid, where $c_1$ and $c_2$ are constants. We also know that the number of unknowns grows linearly within the sequence, in other words $N_l = c_3 + c_4 l$, where $c_3$ and $c_4$ are constants.

### 3.1 NFLOPS estimates for the LU factorization over a sequence of h-refined grids without reutilization

The total NFLOPS required for performing the LU factorization of the entire sequence of *h*-refined grids without reutilizing previously computed LU factorizations is given by:

$$\sum_{l=1}^{L} T_l = \sum_{l=1}^{L} \left( c_1 + N_l c_2 \right) = L c_1 + \sum_{l=1}^{L} N_l c_2 = L c_1 + c_2 \sum_{l=1}^{L} \left( c_3 + c_4 l \right) = L c_1 + L c_2 c_3 + c_2 c_4 \sum_{l=1}^{L} l =$$

$$L c_1 + L c_2 c_3 + c_2 c_4 \left( \frac{L(L+1)}{2} \right) = O\left( L \left( c_1 + c_2 c_3 + c_2 c_4 \right) + L^2 c_2 c_4 \right) = O\left( L + L^2 \right) = O\left( N^2 \right)$$

We can also estimate the memory usage, which is

$$\max_l M_l = \max_l \left( c_5 + N_l c_6 \right) = O(N).$$

### 3.2 NFLOPS estimates for the LU factorization over a sequence of h-refined grids with reutilization

We estimate now the total NFLOPS required for computing the LU factorization over the entire sequence of *h*-refined grids with reutilization of previously computed LU factorizations. When the reutilization is turned on, the cost of solving a single grid from the sequence becomes constant, and we obtain

$$\sum_{l=1}^{L} T_l = \sum_{l=1}^{L} c_1 = L c_1 = O(N).$$

We can also estimate the memory usage. When the reutilization is active, the solution consists of generating just one frontal matrix associated with the root of the elimination tree, which implies that the memory usage of a single grid from the sequence is constant ($M_l = c_7$), and we obtain

$$\sum_{l=1}^{L} M_l = \sum_{l=1}^{L} c_7 = L c_7 = O(N).$$

We conclude that the reutilization reduces the execution time from $O(N^2)$ down to $O(N)$ and preserves the memory usage of order $O(N)$.

## 4. Numerical results

The reutilization solver has been implemented and tested on a number of model problems.
We report here the memory usage expressed as the number of non-zero entries in the LU factorization as well as the execution time for the following solver algorithms:
- the state-of-the-art *MUMPS* direct solver [12,13,14] with *PORD* [21] ordering,
- our in-house solver called *Hypersolver* [24] that employs the elimination trees proposed in this paper but do not reutilize partial LU factorizations, and our in-house solver called *Reutilization* that employs the elimination trees proposed in this paper and also reutilizes previously computed partial LU factorizations.

The reason why we usually employ the MUMPS solver with *PORD* ordering and not the *METIS* [20] ordering is that the PORD ordering provides better performance of the solver on the examples considered in this paper.
These solvers have been tested using the following examples:
- the two dimensional radical mesh with one artificially enforced singularity presented in Figure 1,
- the two dimensional radical mesh with two artificially enforced singularities presented in Figure 4,
- the two dimensional L-shape domain problem described [17, 18, 7], and
- the three dimensional Fichera problem introduced in [8]

In all grids, we employ a polynomial order of approximation $p=5$.

We start by comparing the number of new non-zero entries in the LU factorization for the $h$ refined grid sequence towards point singularites. The comparisons for all model problems are described in Figures 5, 6, 7 and 8. Both MUMPS and our home-made Hypersolver have to recompute the entire LU factorization for each new mesh. Thus, the number of non-zero entries grows within the sequence. However, the Reutilization solver only re-computes non-zero entries of the LU factorization in those newly refined elements surrounding the point singularity.

Figures 9, 10, 11 and 12 compare the execution times of different solvers when applied to our model problems. Our home-made Hypersolver is significantly slower in the pre-asymptotic regime than currently state-of-the-art solvers like MUMPS. This is because our solver has not been optimized. However, the purpose of the paper is to show that the scalability is better, which eventually translates in a better performance even without using optimized codes, as clearly illustrated in the results. From these results, we conclude that the reutilization solver outperforms MUMPS solver, since it delivers constant execution time for each grid from the sequence of grids. In particular, for the most expensive Fichera problem, we also distinguish between the forward elimination and backward substitution parts of the solver, and we conclude that the backward substitution time is negligible here, as expected. However, the computational cost for backward substitution is linear for a single grid of the sequence, but the computational cost constant is relatively small.

Finally, Figures 13, 14, 15 and 16 compare the total execution time for the entire sequence of $h$ refined grids. From the comparison, it follows that only our reutilization solver delivers linear execution time.

## 5. Limitations of the method

The linear computational cost of the direct solver algorithm can be only obtained in presence of point singularities. Thus, we have utilized the 3D Fichera problem to test only a sequence of grids refined towards a central point singularity, as displayed in the left panel of Figure 17. In the case of a sequence of refinements towards point and edge singularities, presented on the right panel of Figure 17, the solver does not exhibit linear computational cost. This is because the number of degrees of freedom added during the edge refinement is no longer constant. The numerical experiments with MUMPS solver suggest that the execution time could be of order $O(N^{1.5})$, see Figure 18.

**Conclusions**
In this paper we presented a new direct solver algorithm that constructs an elimination tree based on the analysis of the structure of the computational mesh. The solver focuses on a class of computational meshes that are *h*-refined towards point singularities. The algorithm delivers linear computational cost of the solution for the entire sequence of meshes that are *h*-refined towards one or several point singularities. This is achieved by using an elimination tree that enables reutilization of previously computed partial LU factorizations over the entire unrefined part of the mesh. The reutilization technique reduces the computational cost of the solution of the entire sequence of *h*-refined grids from $O(N^2)$ down to $O(N)$. These theoretical results are illustrated with several two and three dimensional model problems. Future work includes development of a more general algorithm for construction of the elimination tree that enables reutilization of partial LU factorizations for an arbitrary mesh. We also plan to develop the reutilization algorithm for three dimensional grids with combined point and edge singularities, where we expect a reduction on the computational cost from $O(N^2)$ down to $O(N^{1.5})$.

**Acknowledgement**
The work of the first author was supported by Polish National Science Center grants no NN 519 447739 and  NN 519 405737.
The second author was partially funded by the Project of the Spanish Ministry of Sciences and Innovation MTM2010-16511, the Laboratory of Mathematics (UFI 11/52), and the Ibero-American Project CYTED 2011 (P711RT0278).

**Bibliography**

[1] T. A. Davis, W. W. Hager, Dynamic supernodes in sparse Cholesky update / downdate and triangular solves, ACM Transactions on Mathematical Software, 35, 4 (2009) 1–23.

[2] A. Gupta, G. Karypis, V. Kumar, Highly scalable parallel algorithms for sparse matrix factorization, IEEE Transactions on Parallel and Distributed Systems, 8, 5 (1997) 502-520.

[3] A. Gupta, Recent advances in direct methods for solving unsymmetric sparse systems of linear equations, ACM Transactions on Mathematical Software (TOMS), no. RC 22039 (98933) (2002)

[4] A. Gupta, F. G. Gustavson, M. Joshi, S. Toledo, The design, implementation, and evaluation of a symmetric banded linear solver for distributed-memory parallel computers, ACM Transactions on Mathematical Software, 24,  1 (1998) 74-101.

[5] J. D. Hogg, J. K. Reid, J. A. Scott, A DAG-based Sparse Cholesky Solver for Multicore Architectures RAL-TR-2009-004 (2009).

[6] J. Gaidamour, P. Hénon, A Parallel Direct/Iterative Solver Based on a Schur Complement Approach, 2008 11th IEEE International Conference on Computational Science and Engineering, (2008) 98-105.

[7] L. Demkowicz, Computing with hp-Adaptive Finite Elements, Vol. I. One and Two Dimensional Elliptic and Maxwell Problems. Chapman & Hall / CRC Applied Mathematics & Nonlinear Science, (2006).

[8] L. Demkowicz, J. Kurtz, D. Pardo, M. Paszyński, W. Rachowicz, A. Zdunek, Computing with hp-Adaptive Finite Elements, Vol. II. Frontiers: Three Dimensional Elliptic and Maxwell Problems with Applications, Chapman & Hall / CRC Applied Mathematics & Nonlinear Science, (2007).

[9] Irons B., A frontal solution program for finite-element analysis, International Journal of Numerical Methods in Engineering, 2 (1970) 5-32

[10] I. S. Duff, J. K. Reid, The multifrontal solution of indefinite sparse symmetric linear systems. ACM Transactions on Mathematical Software, 9 (1983) 302-325.

[11] I. S. Duff, J. K. Reid, The multifrontal solution of unsymmetric sets of linear systems. SIAM Journal on Scientific and Statistical Computing, 5 (1984) 633-641.

[12] P. R. Amestoy, I. S. Duff, J.-Y. L'Excellent, Multifrontal parallel distributed symmetric and unsymmetric solvers, Computer Methods in Applied Mechanics and Engineering, 184 (2000) 501-520.

[13] P. R. Amestoy, I. S. Duff, J. Koster, J.-Y. L'Excellent, A fully asynchronous multifrontal solver using distributed dynamic scheduling, SIAM Journal of Matrix Analysis and Applications, 23, 1 (2001) 15-41.

[14] P. R. Amestoy, A. Guermouche, J.-Y. L'Excellent, S. Pralet, Hybrid scheduling for the parallel solution of linear systems, Parallel Computing, 32, 2 (2006) 136-156.

[15] P. Schmitz, L. Ying, A fast direct solver for elliptic problems on general meshes in 2D, Journal of Computational Physics, 231 (2012), 1314-1338.

[16] P. Bientinesi, V. Eijkhout, K. Kim, J. Kurtz, R. van de Geijn, Sparse Direct Factorizations through Unassembled Hyper-Matrices, Computer Methods in Applied Mechanics and Engineering, 199 (2010) 430-438.

[17] I. Babuška, B. Guo, The hp-version of the finite element method, Part I: The basic approximation results, Computational Mechanics, 1 (1986) 21-41.

[18] I. Babuška, B. Guo, The hp-version of the finite element method, Part II: General results and applications, Computational Mechanics, 1 (1986) 203-220.

[19] P. Gurgul, M. Paszyński, A. Szymczak, V. Calo, D. Pardo, A linear complexity direct solver for h-adaptive radical grids, submitted to Computing and Informatics, 2012.

[20] G. Karypis, V. Kumar, MeTis: Unstructured Graph Partitioning and Sparse Matrix Ordering System, University of Minesota, http://www.cs.umn.edu/~metis, 2011.

[21] J. Schulze, Towards a tighter coupling of bottom-up and top-down sparse matrix ordering methods, BIT, 41 (2001) 800-841.

[22] A. Szymczak, M. Paszynski, D. Pardo, Grach grammar based Petri nets controlled direct solver algorithm, Computer Science 11 (2010) 65-79.

[23] M. Paszynski, On the parallelization of self-adaptive hp finite element methods, Part I. Composite Programmable Graph Grammar, Fundamenta Informaticae 93, 4 (2009) 411-434.

[24] M. Paszynski, D. Pardo, A. Paszynska, Parallel multi-frontal solver for p adaptive finite element modeling of multi-physics computational problems, Journal of Computational Science, 1, 1 (2010), 48-54.

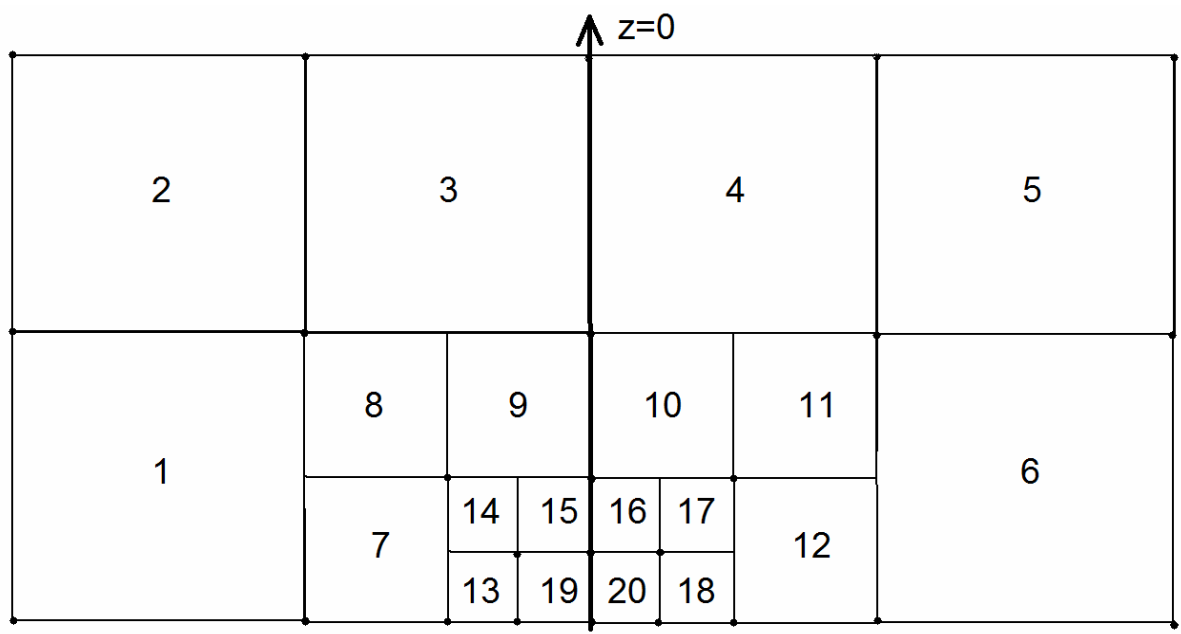

**Figure 1**. Example of a radical mesh.

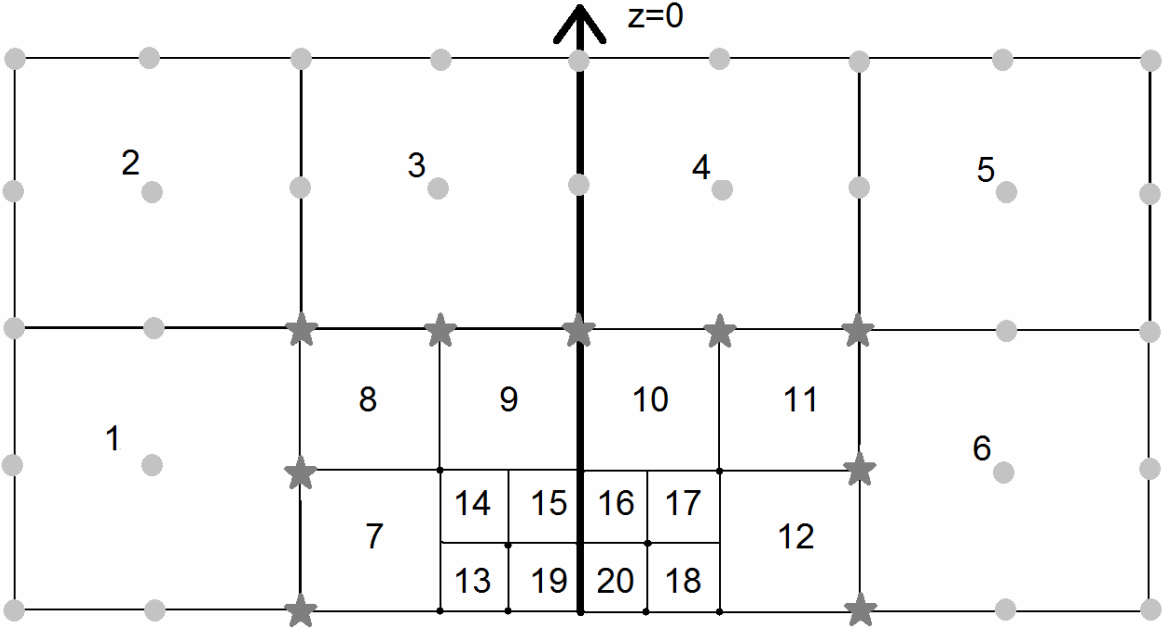

**Figure 2**. Elimination pattern on a single level.

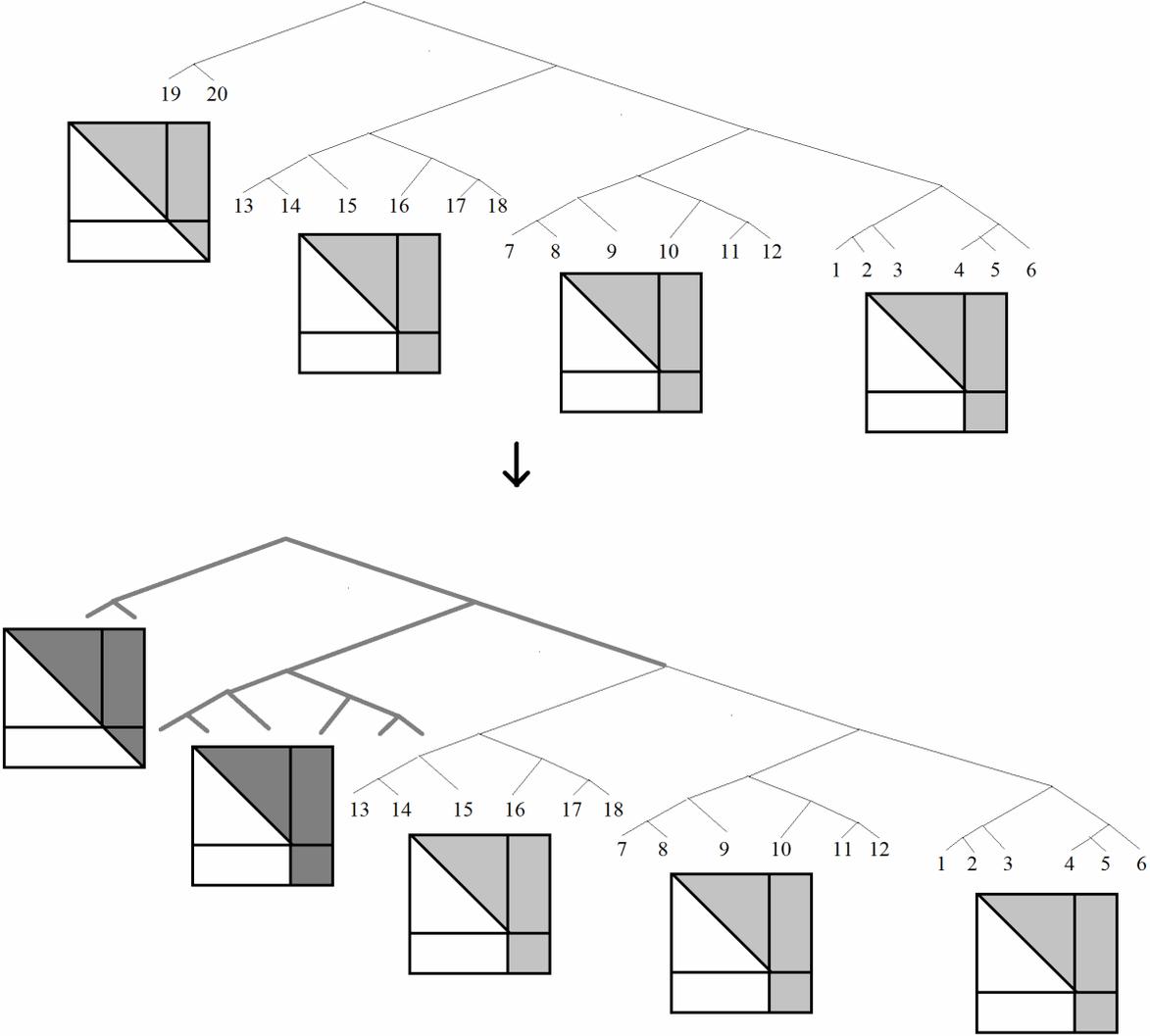

**Figure 3.** Elimination tree based on top-to-bottom ordering, which enables efficient reutilization of previously computed LU factorizations.

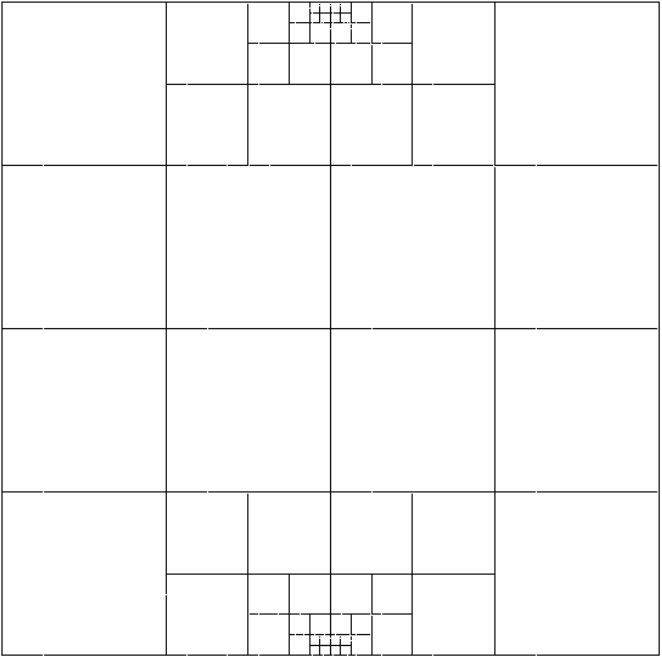

**Figure 4.** Radical mesh with two point singularities.

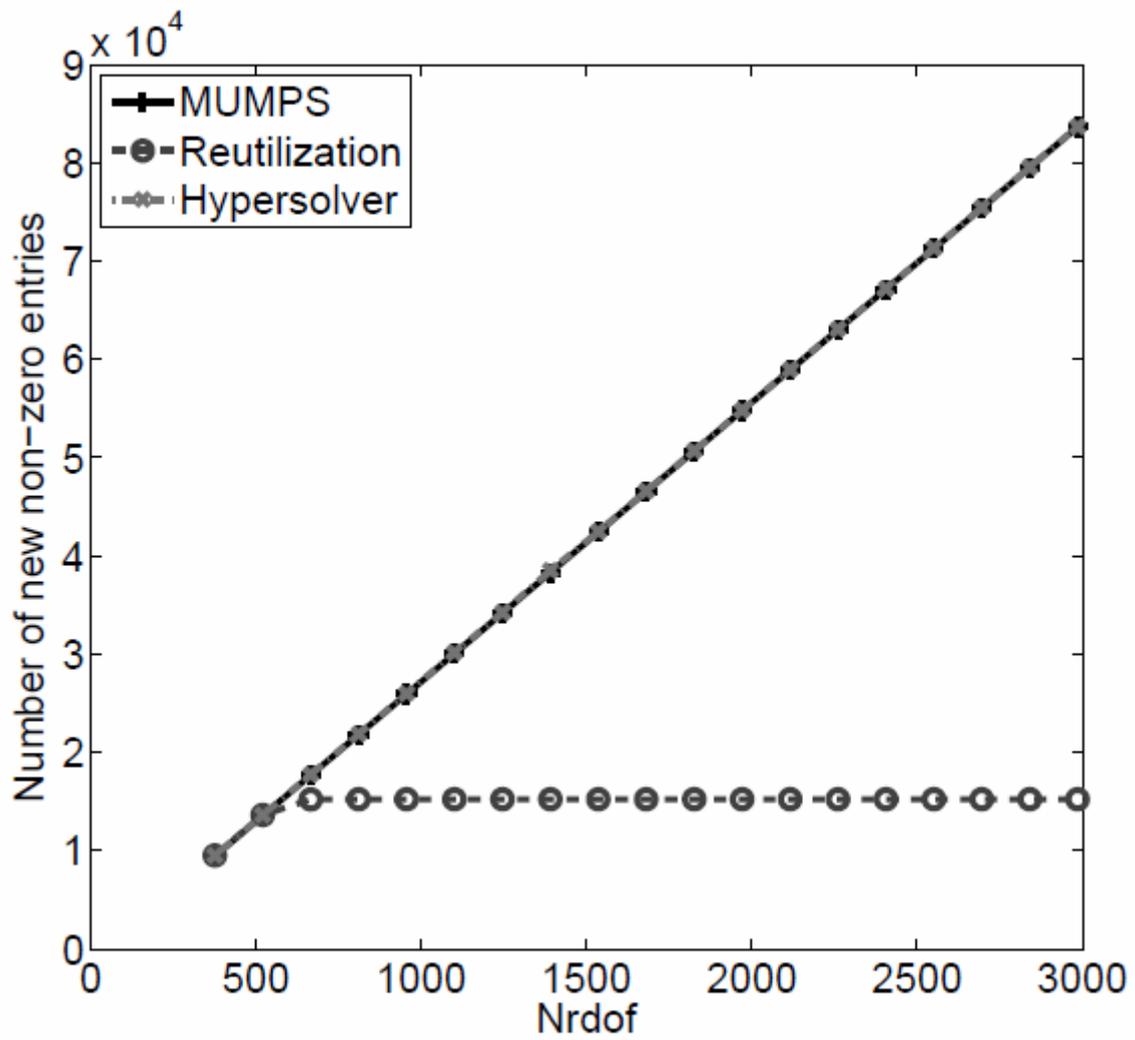

**Figure 5.** Comparison of the number of non-zero entries for the radical mesh with one singularity.

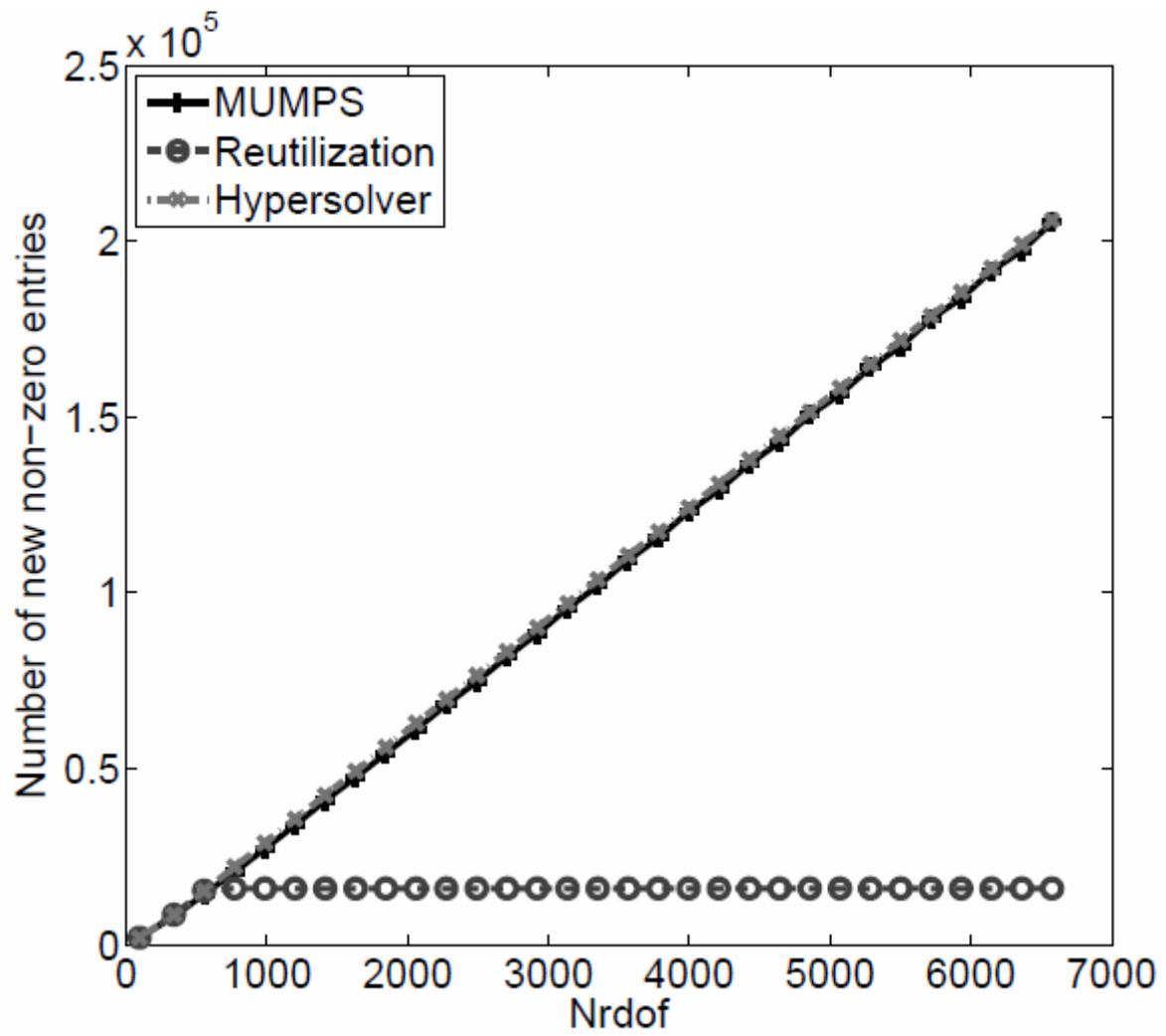

**Figure 6.** Comparison of the number of non-zero entries for the L-shape domain problem.

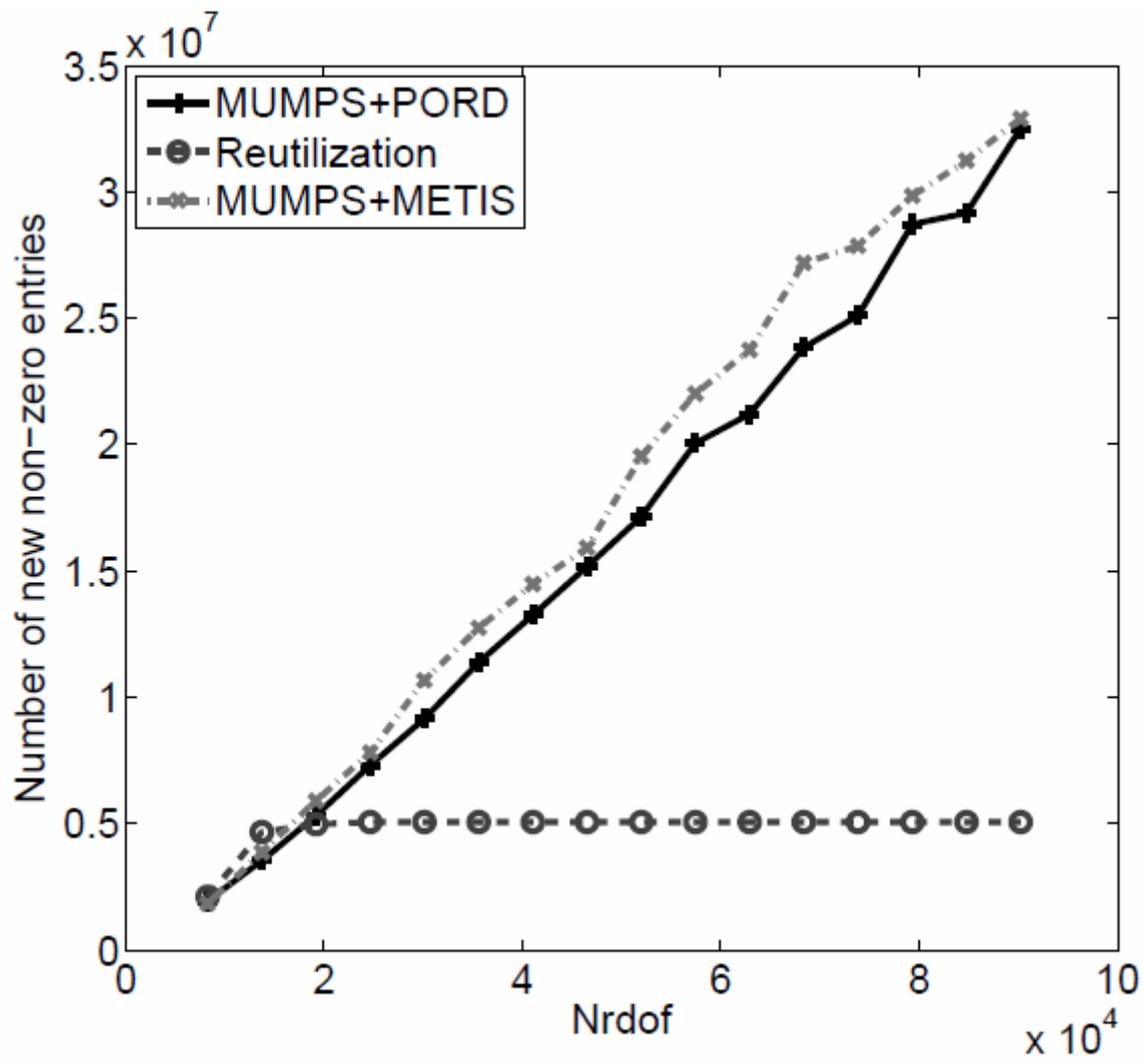

**Figure 7.** Comparison of the number of non-zero entries for the Fichera problem.

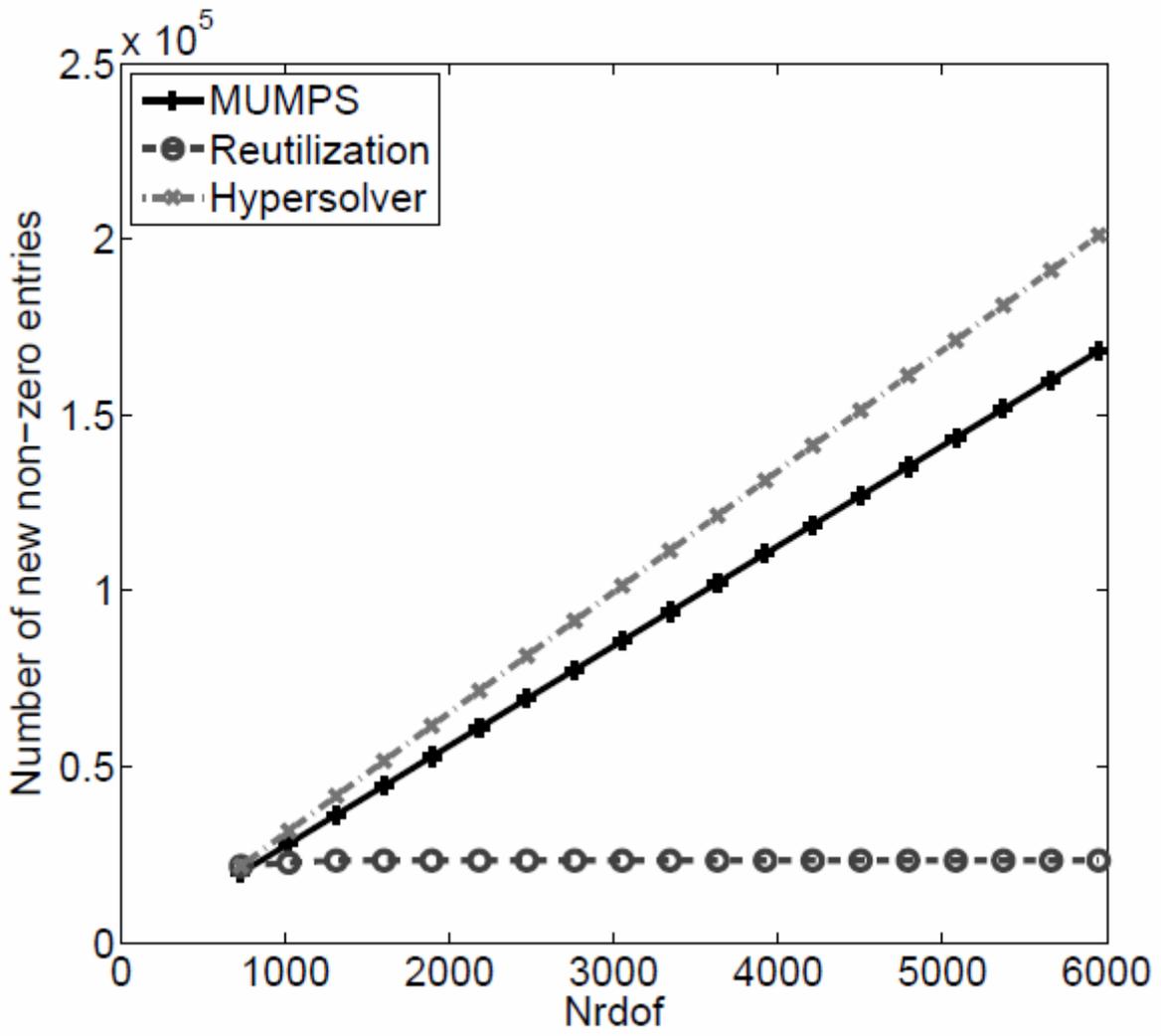

**Figure 8.** Comparison of the number of non-zero entries for the radical mesh with two point singularities.

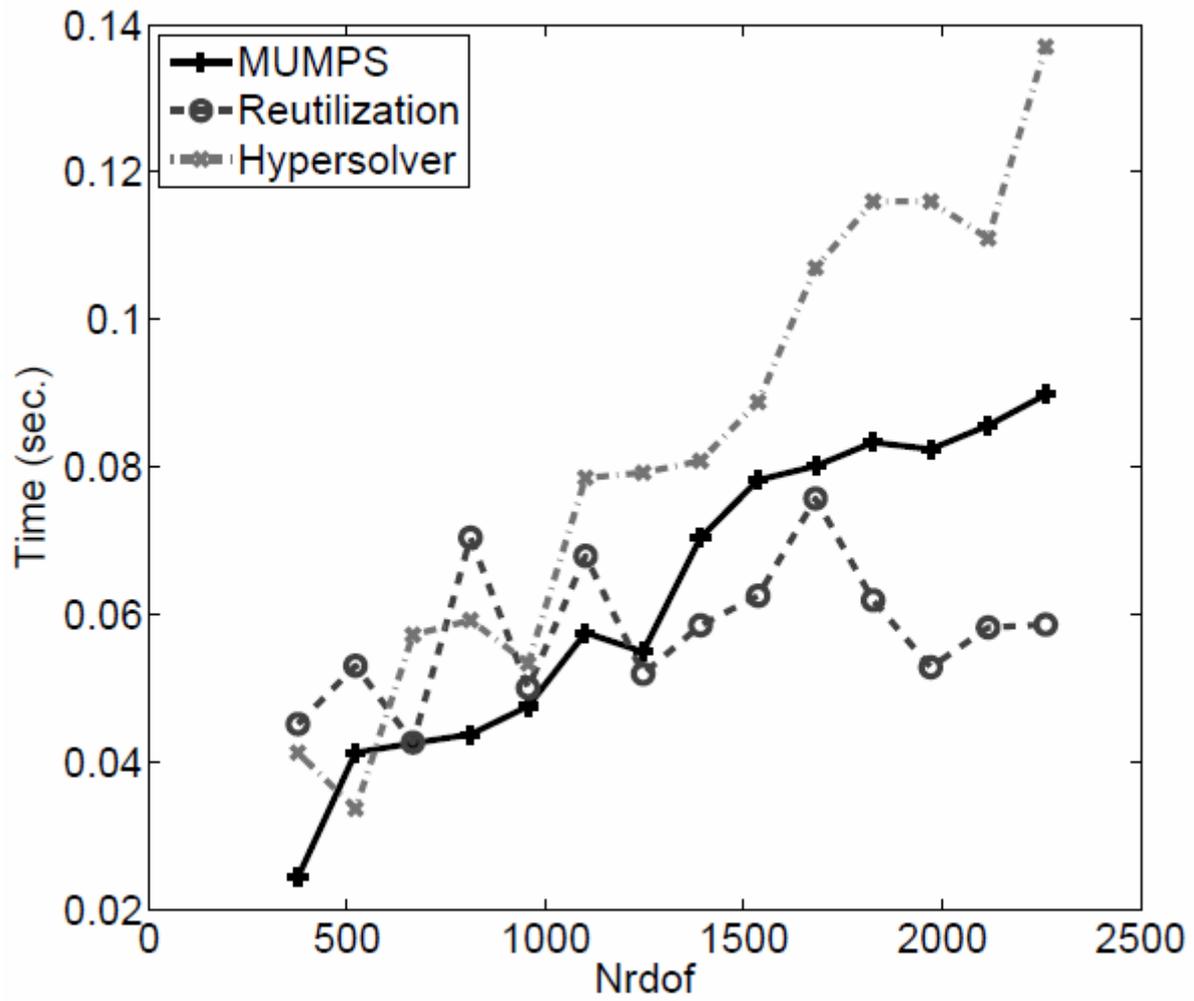

**Figure 9.** Comparison of the execution time for the radical mesh with one point singularity.

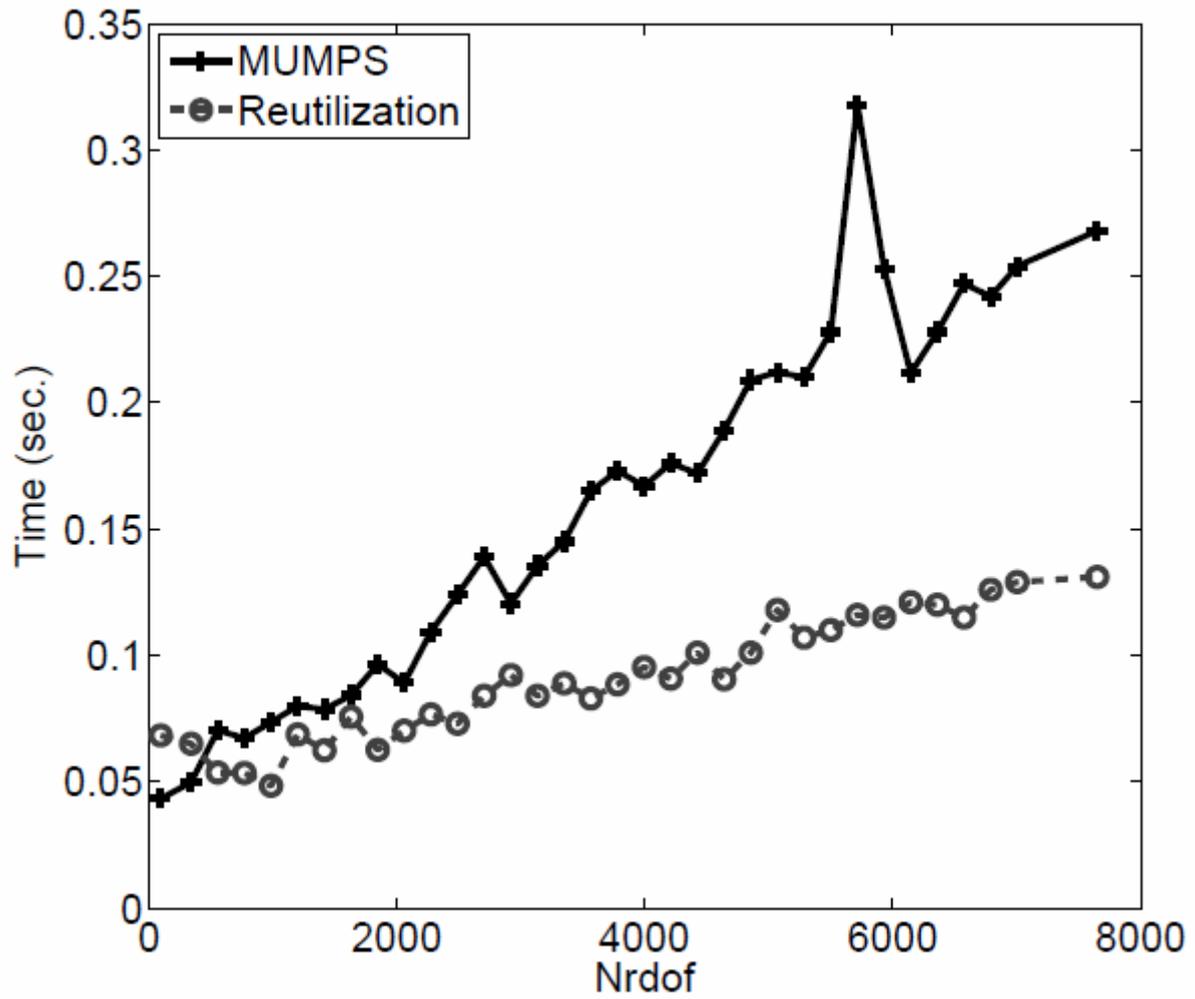

**Figure 10.** Comparison of the execution time for the L-shape domain problem.

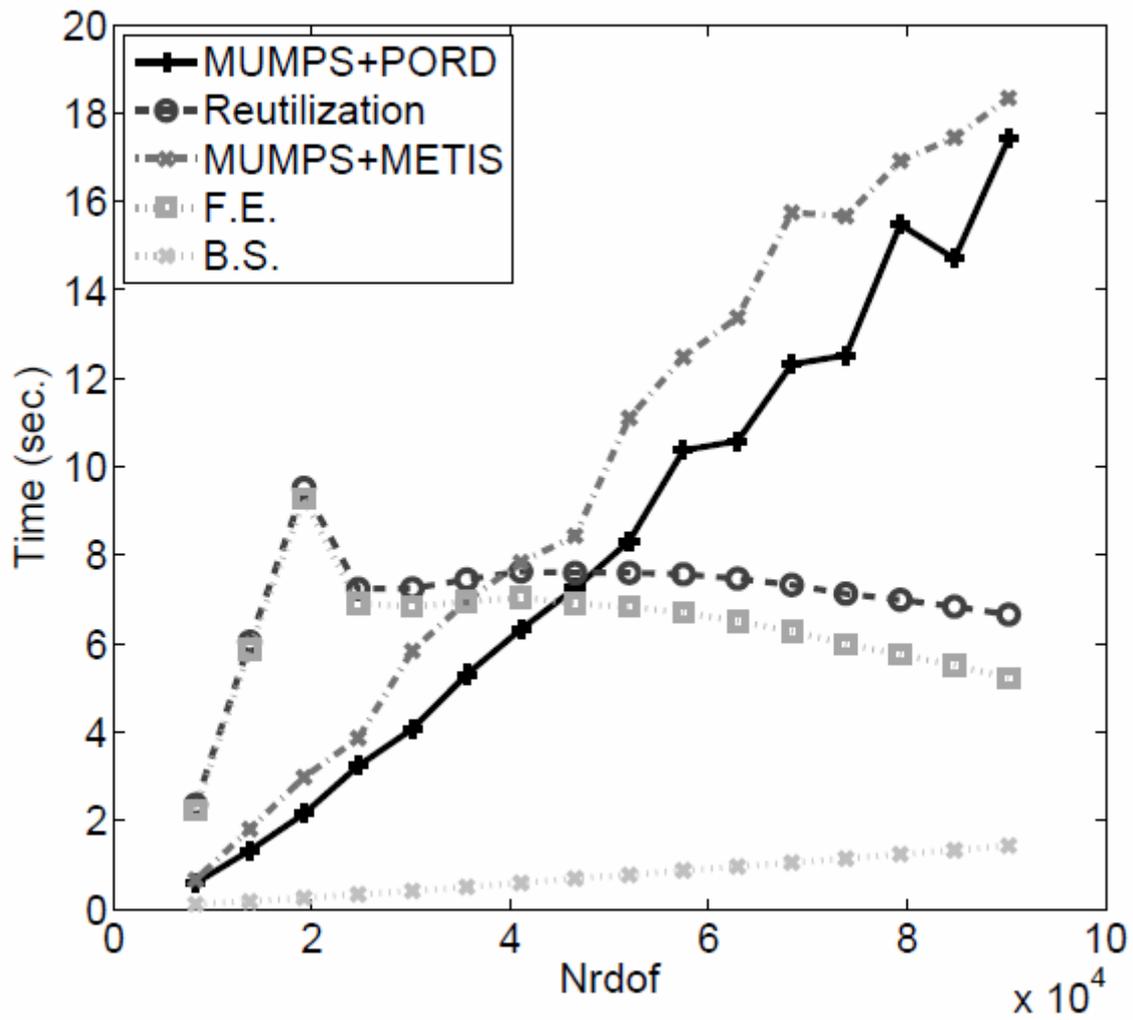

**Figure 11.** Comparison of the execution time for the Fichera problem. Additionally, we describe the forward elimination part (F.E.) and the backward substitution part (B.S.) of the reutilization solver. Notice that the total time of the reutilization solver is the sum of the F.E. and B.S. times.

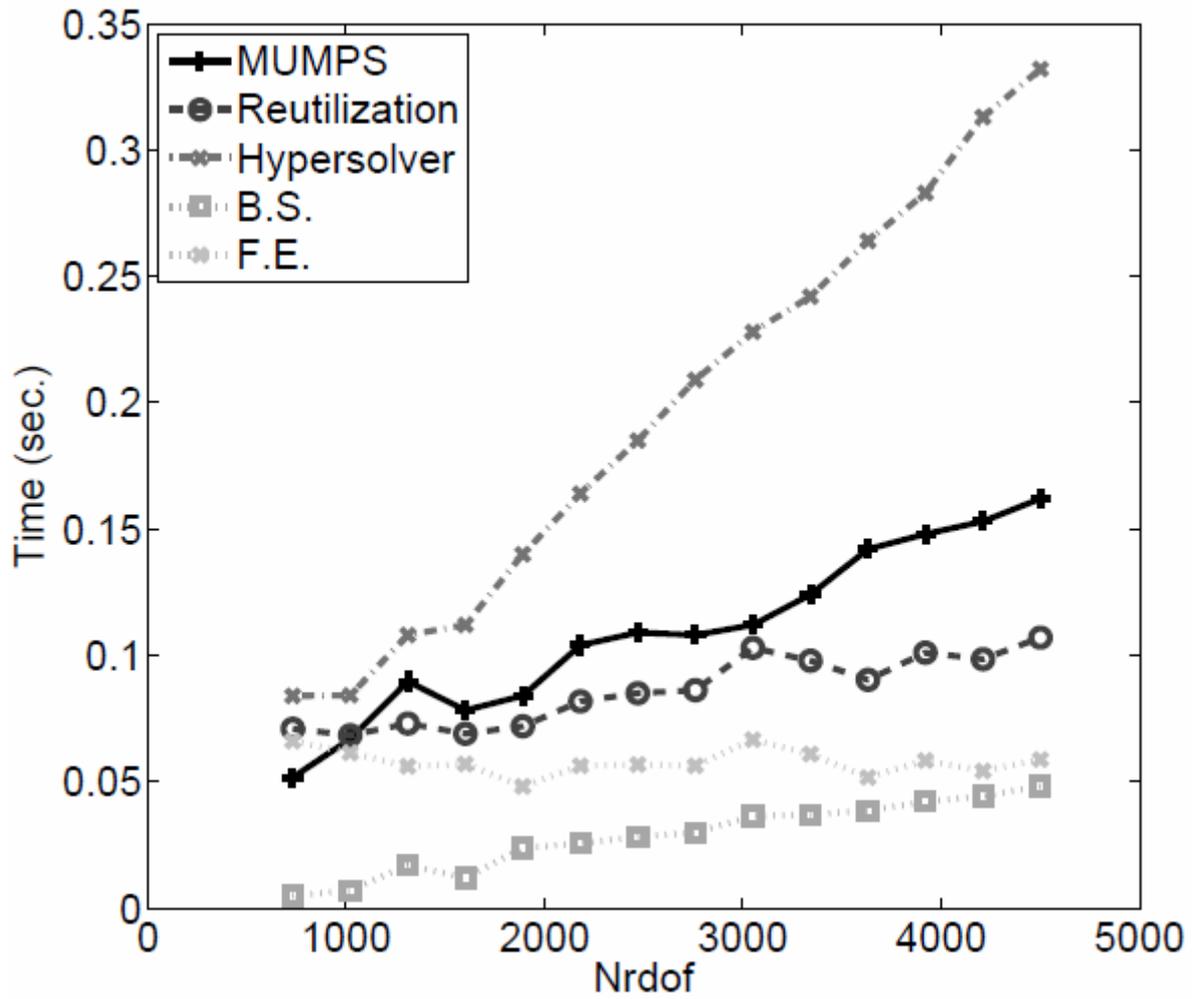

**Figure 12.** Comparison of the execution time for the radical mesh with two singularities. Additionally, we describe the forward elimination part (F.E.) and the backward substitution part (B.S.) of the reutilization solver. Notice that the total time of the reutilization solver is the sum of the F.E. and B.S. times.

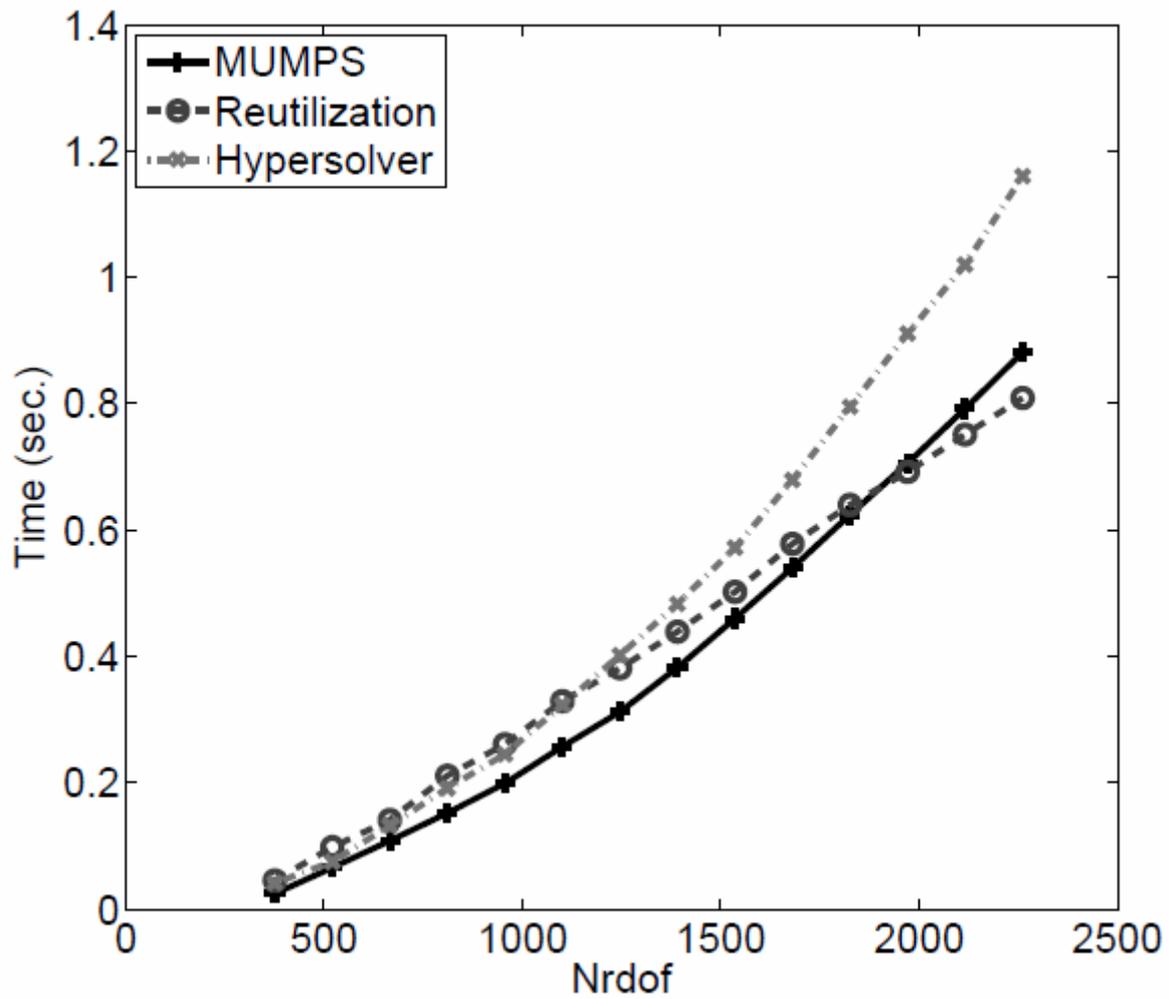

**Figure 13.** Comparison of the total execution time for the entire sequence of grids for the radical mesh with one point singularity.

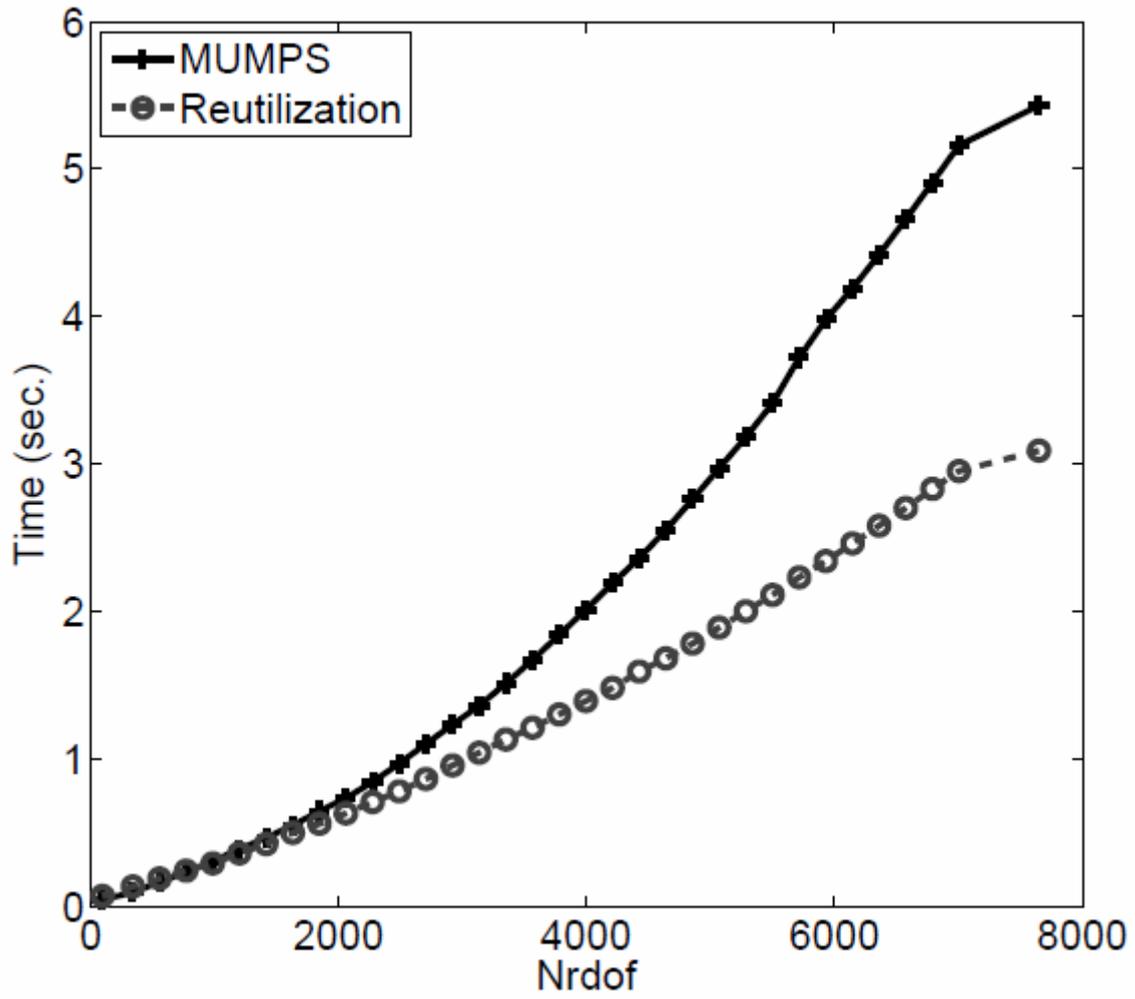

**Figure 14.** Comparison of the total execution time for the entire sequence of grids for the L-shape domain problem.

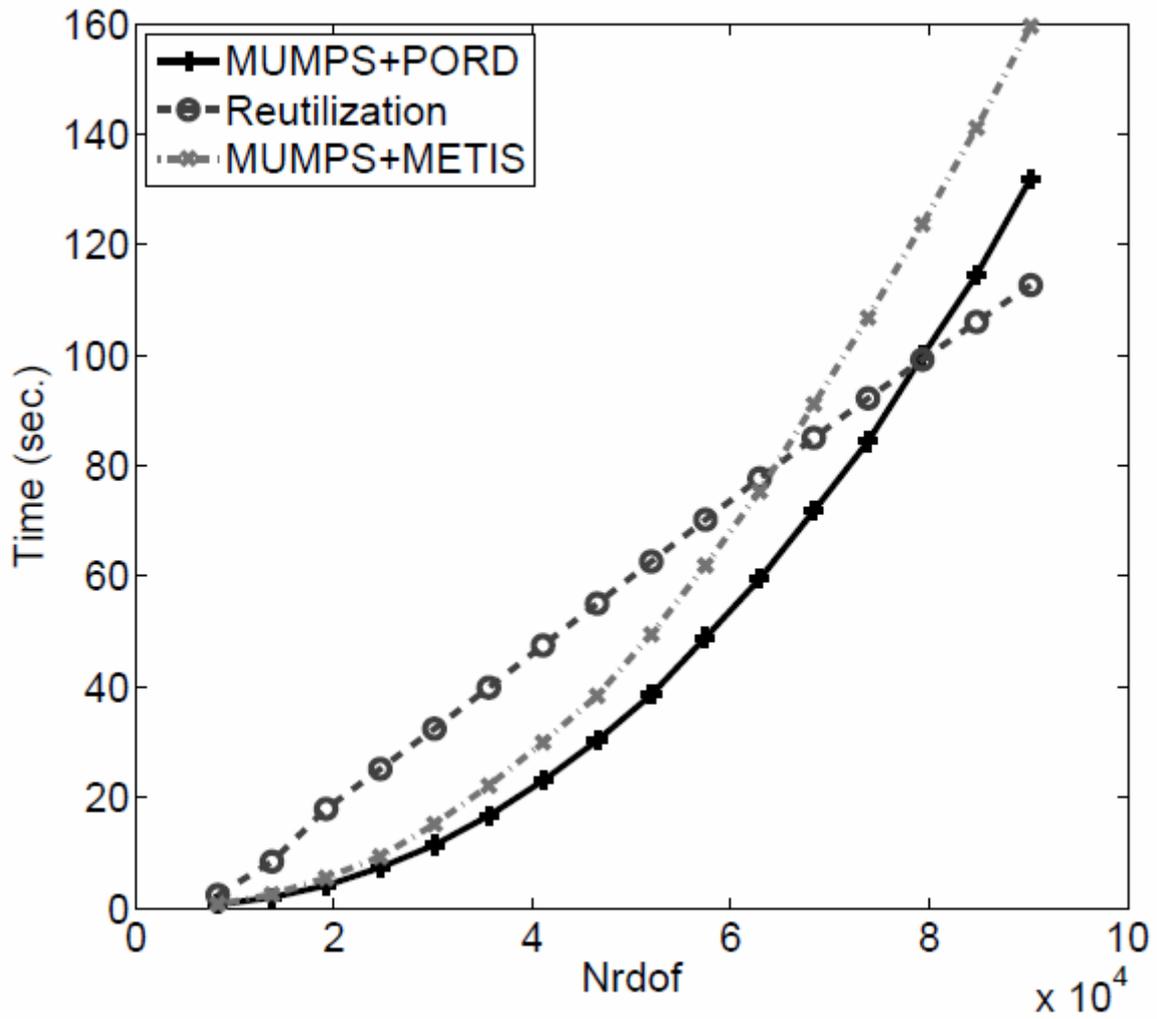

**Figure 15.** Comparison of the total execution time for the entire sequence of grids for the Fichera problem.

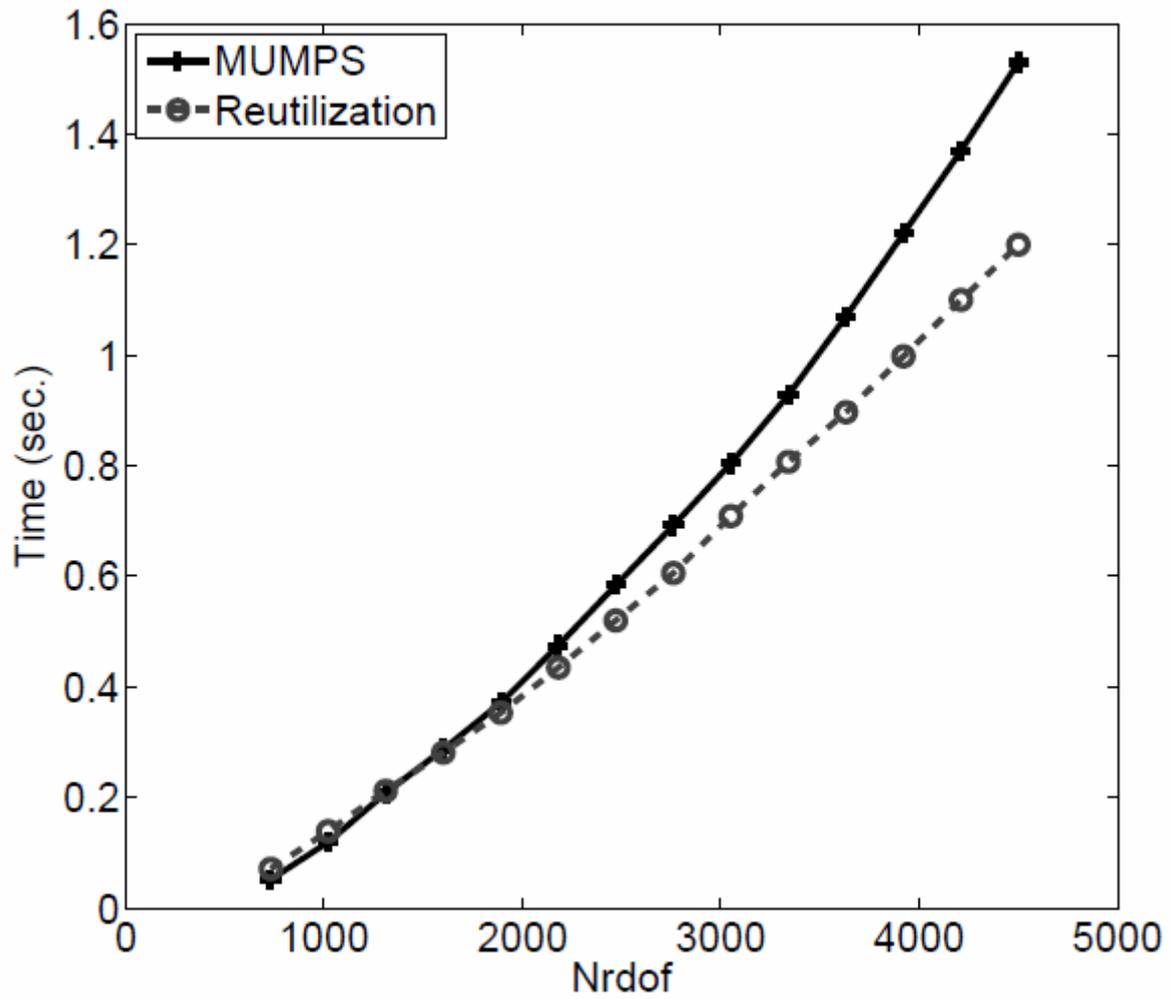

**Figure 16.** Comparison of the total execution time for the entire sequence of grids for the radical mesh with two point singularities.

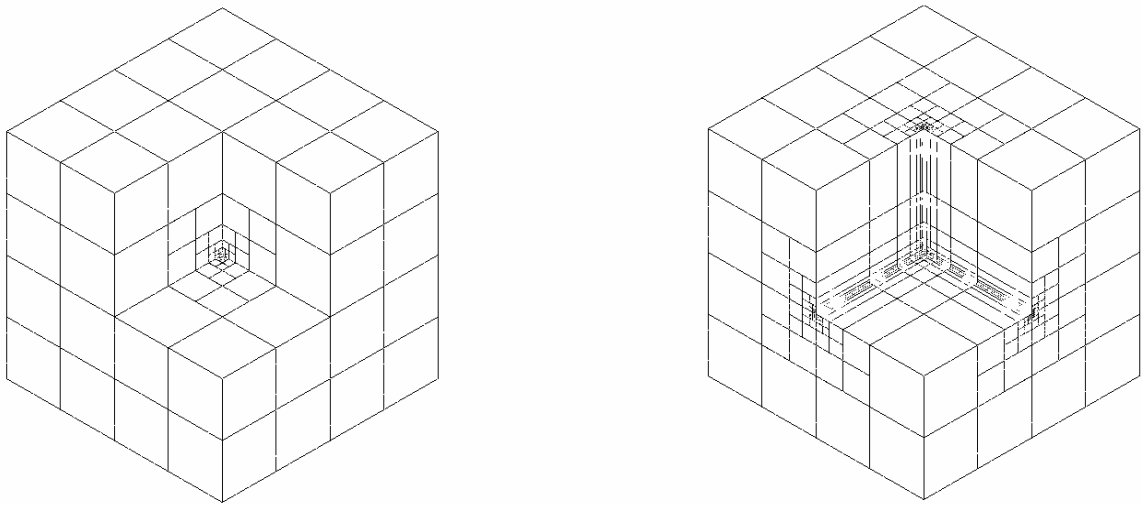

**Figure 17.** Refinements towards point singularity versus refinements towards one point and three edge singularities.

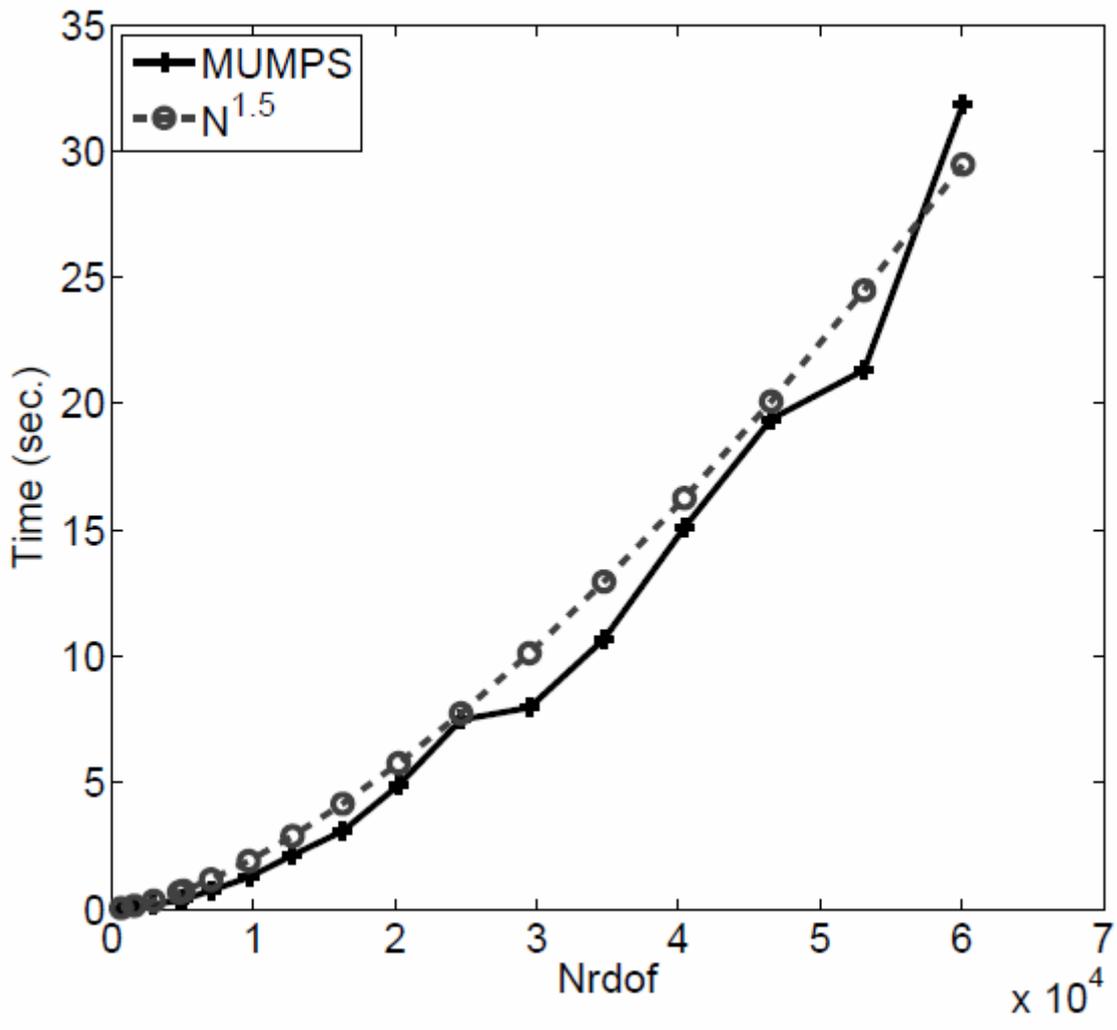

**Figure 18.** Non-linear computational cost of the MUMPS solution corresponding to the Fichera problem with one point and three edge singularities.